\documentclass[12pt]{amsart}
\usepackage{amscd,amssymb}
\begin{document}
\newtheorem{thm}{Theorem}
\newtheorem{lem}{Lemma}
\newtheorem{rem}{Remark}
\newtheorem{conj}{Conjecture}
\newtheorem{prop}{Proposition}
\newtheorem{defn}{Definition}
\newcommand{\ls}{lexicographically shellable}
\newcommand{\lsg}{lexicographical shelling}
\newcommand{\ze}{\mbox{$\hat{0}$}}
\newcommand{\one}{\mbox{$\hat{1}$}}
\newcommand{\hp}{\mbox{$\hat{P}$}}
\newcommand{\bj}{Bj\"{o}rner}
\begin{center}
\large A Note on Planar and Dismantlable Lattices
\end{center}
\begin{center}
Karen L. Collins \\
Dept. of Mathematics, Wesleyan University \\
Middletown, CT 06459-0128 
\end{center}

\begin{abstract}
It is shown that any finite, rank-connected, dismantlable lattice is 
lexicographically shellable (hence Cohen-Macaulay).  
A ranked, interval-connected lattice is
shown to be rank-connected, but a rank-connected lattice need
not be interval-connected.  An example of a planar, rank-connected
lattice that is not admissible is given.
\end{abstract}


\vspace{.03in}

All lattices are assumed to be finite.  \bj\ \cite{B2} has shown that
a dismantlable (see Rival, \cite{R}) lattice $L$ is Cohen-Macaulay 
(see \cite{3} for definition) if and only if $L$ is ranked and 
interval-connected.  
A lattice is planar if its Hasse diagram can be drawn in the plane
with no edges crossing.  Baker, Fishburn and Roberts have shown
that planar lattices are dismantlable, see \cite{BFR}.  
Lexicographically shellable lattices are Cohen-Macaulay, see \cite{B}.
In a recent paper, \cite{C},
the author proved that a planar lattice $L$ is lexicographically shellable
if and only if $L$ is rank-connected.  

\[\begin{picture}(320,115)(0,-35)
\multiput(40,0)(0,20){3}{\circle*{4}}
\multiput(0,20)(0,20){2}{\circle*{4}}
\multiput(20,20)(0,20){2}{\circle*{4}}
\multiput(60,20)(0,20){2}{\circle*{4}}
\multiput(80,20)(0,20){2}{\circle*{4}}
\put(40,80){\circle*{4}}
\multiput(10,60)(20,0){4}{\circle*{4}}
\multiput(0,20)(20,0){5}{\line(0,1){20}}
\multiput(0,40)(20,0){4}{\line(1,2){10}}
\multiput(20,40)(20,0){4}{\line(-1,2){10}}
\multiput(20,20)(20,0){3}{\line(1,1){20}}
\put(40,0){\line(-1,1){20}}
\put(40,0){\line(-2,1){40}}
\put(40,0){\line(1,1){20}}
\put(40,0){\line(2,1){40}}
\put(40,80){\line(-3,-2){30}}
\put(40,80){\line(-1,-2){10}}
\put(40,80){\line(1,-2){10}}
\put(40,80){\line(3,-2){30}}
\put(40,0){\line(0,1){20}}
\put(0,40){\line(4,-1){80}}

\put(140,0){\circle*{4}}
\put(120,20){\circle*{4}}
\put(160,20){\circle*{4}}
\multiput(100,40)(20,0){5}{\circle*{4}}
\multiput(120,60)(20,0){3}{\circle*{4}}
\put(140,80){\circle*{4}}
\put(140,40){\line(0,1){40}}
\put(120,60){\line(1,1){20}}
\put(140,80){\line(1,-1){20}}
\put(100,40){\line(1,1){20}}
\put(120,20){\line(0,1){40}}
\put(160,20){\line(0,1){40}}
\put(120,40){\line(1,1){20}}
\put(140,60){\line(1,-1){20}}
\put(160,60){\line(1,-1){20}}
\put(100,40){\line(1,-1){20}}
\put(120,20){\line(1,1){20}}
\put(160,20){\line(-1,1){20}}
\put(160,20){\line(1,1){20}}
\put(160,20){\line(-1,-1){20}}
\put(120,20){\line(1,-1){20}}

\multiput(200,30)(20,0){3}{\circle*{4}}
\multiput(220,50)(20,0){3}{\circle*{4}}
\multiput(210,10)(20,0){2}{\circle*{4}}
\multiput(230,70)(20,0){2}{\circle*{4}}
\put(220,0){\circle*{4}}
\put(220,0){\line(-1,1){10}}
\put(220,0){\line(1,1){10}}
\put(240,80){\circle*{4}}
\put(240,80){\line(-1,-1){10}}
\put(240,80){\line(1,-1){10}}
\multiput(200,30)(20,0){2}{\line(1,-2){10}}
\multiput(210,10)(20,0){2}{\line(1,2){10}}
\multiput(230,70)(20,0){2}{\line(1,-2){10}}
\multiput(220,50)(20,0){2}{\line(1,2){10}}
\multiput(220,50)(20,0){2}{\line(0,-1){20}}
\multiput(200,30)(40,0){2}{\line(1,1){20}}
\put(220,50){\line(1,-1){20}}

\multiput(280,20)(20,0){3}{\circle*{4}}
\multiput(280,40)(20,0){3}{\circle*{4}}
\put(300,0){\circle*{4}}
\put(300,60){\circle*{4}}
\put(290,40){\circle*{4}}
\put(300,0){\line(0,1){60}}
\multiput(280,20)(40,0){2}{\line(0,1){20}}
\multiput(280,20)(0,20){2}{\line(1,1){20}}
\multiput(300,40)(0,20){2}{\line(1,-1){20}}
\put(300,0){\line(-1,1){20}}
\put(300,0){\line(1,1){20}}
\put(290,40){\line(1,-2){10}}
\put(290,40){\line(1,2){10}}

\put(34,-15){(a)}
\put(134,-15){(b)}
\put(224,-15){(c)}
\put(294,-15){(d)}
\put(145,-35){Figure 1}

\end{picture}\]

We prove a conjecture of \bj\
that a dismantlable, rank-connected lattice is lexicographically
shellable.  We also show that an ranked and interval-connected lattice must be 
rank-connected.  Hence, if $L$ is a dismantlable lattice, $L$ ranked 
and is
interval-connected if and only if $L$ is rank-connected.
However, a rank-connected lattice need not be interval-connected.
Figure 1(a) is a rank-connected lattice that is neither interval-connected 
nor planar.  Not every dismantlable lattice is planar, see, for 
instance, Figure 1(d).

In \cite{C} it was conjectured that planar, rank-connected lattices are
admissible (see Stanley, \cite{S1}).  However, Figure 1(b) is a 
counterexample to that conjecture.  Figure 1(c) is a planar,
rank-connected lattice that is neither upper nor lower semi-modular.

A lattice must have a least element $\ze$, and a greatest element
$\one$.  A lattice that contains only a least element and a greatest
element is trivial.   A lattice is {\bf ranked} if every maximal chain
from $\ze$ to $\one$ has the same length.  
For element $x$, $r(x)$ is defined to be the length of a maximal
chain from the least element to $x$.  Let $R_{i}$ be the set of
elements of rank $i$.  A lattice is {\bf rank-connected} if it is
ranked and the subgraph of the Hasse diagram induced by $R_{i}$ and
$R_{i+1}$ forms a connected graph for all $0\leq i < r(\one)$.  

An element in a lattice is said to be {\bf join-irreducible}
if it covers exactly one element, and {\bf meet-irreducible} if it is
covered by exactly one element.  An element that is both join- and
meet-irreducible is said to be {\bf doubly irreducible}.   
A lattice $L$ is {\bf dismantlable} if there is a chain
$L_{1}\subseteq L_{2} 
\subseteq \cdots \subseteq L_{n}=L$ of sublattices of $L$ such that
the lattice $L_{i}$ has $i$ elements for $1\leq i \leq n$. 
Each single element in $L_{k}-L_{k-1}$ is doubly irreducible
in $L_{k}$.   
Equivalently, a lattice is dismantlable if and only if
every non-trivial sublattice has a doubly irreducible element, see \cite{R}.

Let $L$ be a ranked lattice.  Let $C(L)$ equal the set of covering 
relations of $L$.  Then $L$ is {\bf \ls\ }if there exists a labeling 
$f:C(L)\rightarrow {\bf R}$ such that 
\begin{enumerate}
\item in every interval $[x,y]$ of $L$ there is a unique unrefinable
chain $x=x_{0}<x_{1}<\ldots <x_{n}=y$ such that
$f(x_{0},x_{1})\leq f(x_{1},x_{2})\leq \ldots \leq f(x_{n-1},x_{n})$.
\item for every interval $[x,y]$ of $L$, if $x=x_{0}<x_{1}<\ldots <x_{n}=y$
is the unique unrefinable chain with rising labels, 
and if $z\in [x,y]$ covers $x$ with $z\neq x_{1}$, then $f(x,x_{1})<f(x,z)$.
\end{enumerate}

\begin{thm} \rm Let $L$ be a rank-connected, dismantlable lattice.
Then $L$ is lexicographically shellable.
\end{thm}

\begin{proof} We prove this by induction.  First, we make a reduction.
Suppose that some rank other than the top and bottom rank 
contains a single vertex $x$. 
Then the lattices $[\ze,x],[x,\one]$ are dismantlable and
rank-connected.  Any lexicographic shelling 
of $[\ze,x]$ and separately of $[x,\one]$ will be a lexicographic
shelling of $L$ as long as all labels in $[x,\one]$ are greater
than all labels of $[\ze,x]$.  By induction, therefore, we may assume 
that every rank except the very top and the very bottom 
contains at least two elements.

Since $L$ is dismantlable, every non-trivial sublattice of $L$ contains
a doubly irreducible element.  We show that the covers of doubly 
irreducible elements and the vertices covered by doubly irreducible
elements cannot themselves be doubly irreducible.
Let $x$ be doubly irreducible in $L$.  Suppose $r(x)=i$.  
Let the unique vertex that $x$ covers be $z$ and the unique vertex
that covers $x$ be $y$.
Let $G$ be the induced bipartite graph of the Hasse
diagram with vertices $R_{i}\cup R_{i+1}$.  Then both $x$ and $y$
are in $G$, and $G$ is a connected graph because $L$ is rank-connected.  
Since $x$ is doubly irreducible,
$x$ is a degree 1 vertex in $G$, hence $G-x$ is still connected.
Since $R_{i}$ contains at least two
vertices, $G-x$ must contain $y$ and at least one vertex in $R_{i}$,
hence $y$ has an edge to another vertex of rank $i$.  Therefore
$y$ (and similarly $z$) cannot be doubly irreducible in $L$.

For rest of the proof, we rely on the following definition and previous result.  
We say that $w$ is a corner of $x$ in $L$ if there exist $z$ and $y$ such 
that $x,w$ both cover $z$ and are covered by $y$, and $x$ is  
doubly irreducible.  Theorem~1 of \cite{C} 
guarantees that $L$ is lexicographically shellable
if $L-x$ is lexicographically shellable.   
We prove that there must always exist a 
doubly irreducible element of $L$ that has a corner, hence
$L$ is lexicographically shellable by induction.

Let the doubly irreducible vertices of $L$ be $D:=\{x_{1},x_{2},x_{3},
\ldots,x_{t}\}$.  Suppose that no element of $D$ has a corner.  
Let $x_{j}$ cover $z_{j}$ and be covered by $y_{j}$ in $L$ for
$1\leq j\leq t$.  As seen before, rank connectedness guarantees that
$z_{j}$ and $y_{j}$ are not doubly irreducible in $L$.  
Consider the relation $z_{j}\leq y_{j}$ (which is true in $L$) in the
sublattice $L-D$.  There must be a chain from $z_{j}$ to $y_{j}$
in $L-D$.

However in $L$, the rank of $z_{j}$ and $y_{j}$ differs by exactly
2.  Therefore, any chain between $z_{j}$ and $y_{j}$ in $L-D$
must have length less than or equal to 2, since removing vertices
cannot make chains longer.  If the length of the chain is 2, then
the middle vertex of the chain will be a corner of $x_{j}$.  By our
assumption, $x_{j}$ has no corners, hence $y_{j}$ must cover $z_{j}$ 
in $L-D$.

Suppose that $y_{j}$ covers both $x_{j}$ and $x_{k}$ for some $k$.
Now $x_{k}$ is not a corner of $x_{j}$, so $z_{j}\neq z_{k}$.  Therefore,
$y_{j}$ covers both $z_{j}$ and $z_{k}$ in $L-D$.
Thus if element $u$ covers $s$ vertices in $L$, after removing
the doubly irreducible vertices, $u$ still covers $s$ vertices
in $L-D$.  Each doubly irreducible vertex $x_{j}$ is replaced by 
the unique vertex $z_{j}$ that it covers.  
Similarly, if element $u$ is covered by $s$ vertices in $L$, 
it is still covered by $s$ elements in $L-D$. 

Hence the sublattice $L-D$ is composed entirely of elements
that are not doubly irreducible.  It is not empty, since $L$ must
contain at least one doubly irreducible element $x$ and therefore
contains both the unique element that covers and the unique 
element that is covered by $x$.  Since $L$ is dismantlable, 
$L-D$ must be the trivial lattice that contains only a top
and bottom element.  Thus, $L$ must have rank 2 and all its doubly
irreducible elements have rank 1.  Each of these rank 1 elements is
a corner of every other, contradicting our assumption that no doubly
irreducible element has a corner.  
\end{proof}

Define $[x,y]=\{z|x\leq z\leq y\}$.  A lattice is {\bf interval-connected} if
for every pair $x,y$ with $r(y)\geq 2+r(x)$ the Hasse diagram of
$[x,y]-\{x,y\}$ is connected.  

\begin{thm} \rm Let $L$ be an interval-connected, ranked lattice.  Then
$L$ is rank-connected.
\end{thm}

\begin{proof} Any sublattice $[x,y]$ of $L$ must be interval-connected.
Assume by induction that every interval $[x,y]\neq 
[\ze,\one]$ is rank-connected.  We will show that
$[\ze,\one]$ is rank-connected.  

Let $\hat{L}$ be the subposet of $L$ that contains all vertices
except $\ze, \one$.  Let $G$ be the
subgraph of the Hasse diagram induced by $R_{i}$ and $R_{i+1}$.
Since $L$ is interval-connected, there must be some path between any
two vertices in $G$ in $\hat{L}$. 
Suppose that $G$ is not connected.  Let $u,v$ be elements of $G$ 
such that $u$ and $v$ are not in the same connected component of $G$
and such that $p(u,v)$ is the shortest possible path between any
two elements in different connected components of $G$ that is 
contained in $\hat{L}$.  Let $P$ be the set of vertices in $p(u,v)-\{u,v\}$.
Then $P$ does not intersect the vertices in 
$G$, so $u,v$ both have rank $i$ or both have rank $i+1$, and the rank 
of all vertices in $P$ must be less than $i$ or greater than $i+1$,
respectively.  Without loss of generality, assume that $u,v$ have rank
$i$ and let $a$ be the last element of $p(u,v)$ that is strictly
less than $u$.

We observe that if $u\wedge v>\ze$, then $[u\wedge v,\one]$ is smaller than $L$,
and therefore is rank-connected by induction.  This means there is a
path from $u$ to $v$ in $G\cap [u\wedge v,\one]$, which is
clearly a path in $G$.  Therefore we may assume that $u\wedge
v=\ze$.  Similarly, we can assume $u\vee v=\one$.  
Now $a<u$ and $a\neq v$ since $r(v)=r(u)$.  Let $b$ be the next element after 
$a$ in $p(u,v)$.  Then
$a$ and $b$ are comparable, hence $b\geq a$, since $u \not \geq b$. If
$v\geq b\geq a$, then $\ze = u\wedge v\geq a >\ze$, a contradiction.
So $b$ is incomparable to both $u$ and $v$.  

There exists a path from $b$ to $\one$ in $L$ that is strictly
rank increasing.  Let $p(b,\one)$ be such a path.  
Now $r(b)=r(a)+1\leq r(u)=i$.  Let $b(i)$ be the element on
$p(b,\hat{1})$ that has rank equal to $i$.  Then $b(i)\geq b\geq
a$ and $u\geq a$ imply that $b(i)\wedge u\geq a>\hat{0}$.  Therefore,
$b(i)$ and $u$ must be in the same connected component of $G$.  

If $b(i)$ is in the same connected component of
$G$ as $v$, then by transitivity $u$ and $v$ are in the same
component.  If $b(i)$ and $v$ are in different components, then we
replace the pair $u$ and $v$ by the pair $b(i)$ and $v$ and take a
path from $b(i)$ to $v$ that consists of
starting at $b(i)$ and following a strictly rank decreasing path to
$b$ and then following the portion of $p(u,v)$ from $b$ to $v$.  This
must be a shorter path than $p(u,v)$ between a rank $i$ element and a rank
$i+1$ element in different connected components of $G$, because
the portion of $p(u,v)$ from $u$ to $b$ goes through $a$ and is 
therefore longer than the strictly rank decreasing path from $b(i)$ 
to $b$.  We bypass $a$, and $a$ has lower rank than $b$.  This contradicts 
the selection of $p(u,v)$ as the shortest possible path. 
\end{proof}

Admissible lattices are \ls, but not all \ls\ lattices are admissible.
See Stanley's paper for details \cite{S1}.  We define admissible
lattices and show that there is a planar, rank-connected lattice which
is lexicographically shellable, but not admissible.
Let $J$ be the set of join-irreducibles of a
lattice.   Define a natural labeling $\omega$ of $J$ to be a map 
$\omega:J \rightarrow {\bf N}$ where ${\bf N}$ is the positive
integers such that if $z,w\in J$ and $z\leq w$, then $\omega(z)\leq
\omega(w)$.  Let $\gamma$ be derived
from $\omega$ by $\gamma(x<y)=min\{\omega(z) | z\in J,x<x\vee z
\leq y\}$.  A lattice $L$ is {\bf admissible} 
if whenever $x<y$ in $L$, there is a unique unrefinable chain
$x=x_{0}<x_{1}<\ldots <x_{m}=y$ such that $\gamma(x_{0},x_{1})\leq
\gamma(x_{1},x_{2})\leq \ldots \leq \gamma(x_{m-1},x_{m})$.

The planar, rank-connected lattice in Figure 1(b) is lexicographically
shellable, but it is not admissible.  Label the vertices 
with $0$ to $11$ starting at the lowest rank and
moving left to right.  Then the join irreducibles are $1,2,3,4,6,7$.
Suppose we have an admissible labeling of the poset.
Then $\omega(3),\omega(4)\geq \omega(1)$, hence $\gamma(10<11)=\omega(1)$.
Clearly $\gamma(7<10)=\omega(6)$.  Since there is a unique chain from
$7$ to $11$, it must be rising and $\omega(6)\leq \omega(1)\leq
\omega(4)$.  By left-right symmetry,  $\omega(4)\leq \omega(2)\leq
\omega(6)$.  Thus $\omega(1)=\omega(2)$, which cannot happen, since
the two chains from $0$ to $5$ will then both be rising. 

The author would like to thank the referees for their
helpful suggestions.

\end{document}